\documentclass[11pt]{article}

\usepackage{a4wide}
\usepackage{amsfonts} 
\usepackage{amsmath}
\usepackage{amssymb}
\usepackage{mathrsfs}
\newtheorem {Def}{Definition}[section]
\newtheorem {Thm}[Def]{Theorem}%[section]
\newtheorem {Prop}[Def]{Proposition}%[section]
\newcommand{\CC}{{\mathbb C}}
\newcommand{\id}{\mathrm{id}}

\newcommand{\ph}{\varphi}
\newcommand{\tens}{\otimes}
\newcommand{\Sum}{\sum\limits}
\newcommand{\M}[1]{M\!\left(#1\right)}
\newcommand{\Span}[1]{\mathrm{span}\left\{#1\right\}}
\newcommand{\refeq}[1]{{\rm (\ref{#1})}}
\newcommand{\proof}{\noindent{\sc Proof.~}}
\newcommand{\qed}{\hfill{\sc Q.E.D.}\medskip}
\newcommand{\Aa}{{\mathscr A}}
\newcommand{\Bb}{{\mathscr B}}
\newcommand{\conv}{\star}
\newcommand{\Aah}{{\Aa\sp{\:\,\!\!\widehat{}}}}
\newcommand{\sh}{{\,\sharp}}
\newcommand{\st}{{\,*}}
\newcommand{\del}{\delta}

\begin{document}

\title{A note on multipliers of discrete quantum groups}
\author{Piotr Miko\l{}aj So\l{}tan\thanks{
Partially supported by Komitet Bada\'n Naukowych grant No.~2PO3A04022, the
Foundation for Polish Science and Deutsche Forschungsgemeinschaft.} \\
\small Department of Mathematical Methods in Physics,\\
\small Faculty of Physics, University of Warsaw\\
\small\tt piotr.soltan@fuw.edu.pl}
\date{}
\maketitle

\abstract{We investigate the problem whether a given multiplier of a tensor 
product of two algebras belongs to the tensor product of multiplier algebras.
We give a characterization of such multipliers in the case when one of the 
algebras is the algebra of functions on a discrete quantum group.}

\section{Introduction}

Let $\Aa$ be an algebra (over $\CC$) with non degenerate product. Then we can
define the multiplier algebra $\M{\Aa}$ as the vector space of pairs
$(\lambda,\rho)$ of linear maps $\Aa\to\Aa$ such that 
\[
\rho(a)b=a\lambda(b)
\]
for all $a,b\in\Aa$. It is customary to treat $\lambda$ and $\rho$ as left and
right multiplication by an auxiliary object $m$, i.e.~$\lambda(a)=ma$ and
$\rho(a)=am$ for all $a\in\Aa$. The element $m$ is called a {\em multiplier}\/
of $\Aa$, while the pair $(\lambda,\rho)$ is traditionally referred to as a
{\em double centralizer}\/ of $\Aa$. Usually the a multiplier $m$ and the map
$\lambda$ corresponding to left multiplication by $m$ are denoted by the same
symbol. Note that if $\Aa$ has a unit then
$\M{\Aa}=\Aa$. 

If $\Bb$ is another algebra with non degenerate product then $\Bb\tens\Aa$ is 
again an algebra with non degenerate product (\cite[Lemma A.2]{mha})
and we can define the multiplier algebra $\M{\Bb\tens\Aa}$. The tensor product
of multiplier algebras of $\Bb$ and $\Aa$ embeds naturally into
$\M{\Bb\tens\Aa}$ (\cite[Prop.~A.3]{mha}). In \cite{ap} the image of this 
embedding was used to characterize almost periodic elements for a discrete 
quantum group and helped in constructing the analogue of Bohr compactification
for discrete quantum groups. Our aim in the present paper is to provide a 
criterion characterizing elements of the image of this embedding.  

For any vector space $V$ we shall denote by $V^\sh$ the space of all linear 
functionals on $V$. If $\Aa$ is an algebra then $\Aa^\sh$ is an $\Aa$-bimodule 
in a natural way: for $f\in\Aa^\sh$ and $a\in\Aa$
\[
\begin{array}{r@{\;=\;}l@{\smallskip}}
\bigl(af\bigr)(b)&f(ba),\\
\bigl(fa\bigr)(b)&f(ab)
\end{array}
\]
for all $b\in\Aa$. An important subspace of $\Aa^\sh$ is the space of
{\em reduced functionals}\/ on $\Aa$ which is by definition
\[
\Aa^\st=\Span{\bigl.afb:\:a,b\in\Aa,\:f\in\Aa^\sh\bigr.}.
\]
Any functional in $\Aa^\st$ admits a natural extension to $\M{\Aa}$. If $\Bb$ 
is another algebra with non
degenerate product then for all $\zeta\in\Bb^*$ and $\xi\in\Aa^\st$ the tensor
product $\zeta\tens\xi$ is a reduced functional on $\Bb\tens\Aa$ and so it
extends to $\M{\Bb\tens\Aa}$. We have the following simple result (\cite{ap}):

\begin{Prop}\label{sli}
Let $\Aa$ and $\Bb$ be algebras with non degenerate products and let $Y$ be a
multiplier of $\Bb\tens\Aa$. Then for any $\xi\in\Aa^\st$ there exists a unique 
multiplier $m\in\M{\Bb}$ such that
\[
(\zeta\tens\xi)(Y)=\zeta(m)
\]
for all $\zeta\in\Bb^\st$.
\end{Prop}

The multiplier $m$ constructed in Proposition \ref{sli} is called a\/ {\em right
slice of $Y$ with $\xi$} and will be denoted by $(\id\tens\xi)(Y)$.

Multiplier algebras have been introduced into the theory of Hopf algebras by 
Van Daele in \cite{mha} where the notion of a multiplier Hopf algebra was
defined. The theory of multiplier Hopf algebras provides a natural framework to
study a variety of quantum groups. In particular discrete quantum groups have
been studied in this framework in \cite{dqg}. 

A multiplier Hopf algebra is a pair $(\Aa,\del)$ of an algebra with non
degenerate product and a homomorphism $\del\colon\Aa\to\M{\Aa\tens\Aa}$ such
that the maps
\[
\begin{array}{r@{\;\colon\Aa\tens\Aa\ni(a\tens b)\longmapsto\;}l@{\smallskip}}
T_1&\del(a)(I\tens b),\\
T_2&(I\tens a)\del(b)
\end{array}
\]
have image equal to $\Aa\tens\Aa$, are bijective onto this image and
$T_2\tens\id$ and $\id\tens T_1$ on $\Aa\tens\Aa\tens\Aa$ commute. 
By composing $\del$ with an
extension of the flip map ($a\tens b\mapsto b\tens a)$ to $\M{\Aa\tens\Aa}$ we
obtain another homomorphism $\del'$. If $(\Aa,\del')$ is a multiplier Hopf
algebra then the multiplier Hopf algebra $(\Aa,\del)$ is called {\em regular}.
If $\Aa$ is a $*$-algebra and $\del$ is a $*$-homomorphism (in this case
$\M{\Aa}$ carries a natural involution) then we say that $(\Aa,\del)$ is a
multiplier Hopf $*$-algebra. 

A {\em discrete quantum group}\/ is a multiplier Hopf $*$-algebra
$(\Aa,\del)$ such that $\Aa$ is a direct sum of a family of full matrix
algebras. Since multiplier Hopf $*$-algebras are automatically regular
(\cite[Sect.~5]{mha}) we see that discrete quantum groups are regular
multiplier Hopf algebras. We shall not make use of the involutive structure of
discrete quantum groups, but regularity will be of importance for our results. 

Discrete quantum groups appeared first in the $\mathrm{C}^*$-algebraic context 
in \cite{pw}.

Let $(\Aa,\del)$ be a multiplier Hopf algebra. A functional $\ph$ on $\Aa$ is
called {\em left invariant}\/ if the map $\Aa\ni
a\mapsto(\id\tens\ph)\del(a)\in\M{\Aa}$ satisfies
\[
(\id\tens\ph)\del(a)=\ph(a)I.
\]
If a non trivial left invariant functional exists then it is unique up to 
multiplication by a constant. Moreover if a left invariant functional exists
then so does a right invariant functional (defined similarly with an obvious
modification). Discrete quantum groups have invariant functionals (\cite{dqg})
and we will freely use the results of the theory of regular multiplier Hopf
algebras with invariant functionals developed in \cite{afgd}. 

\section{Reduced functionals as multipliers}

Let $(\Aa,\del)$ be a multiplier Hopf algebra. It was proved in
\cite[Prop.~6.2]{mha} that the comultiplication $\del$ defines an associative 
algebra structure on $\Aa^\st$ by
\[
(\xi_1\conv\xi_2)(a)=(\xi_1\tens\xi_2)\bigl(\del(a)\bigr)
\] 
for all $\xi_1,\xi_2\in\Aa^\st$ and all $a\in\Aa$. We call this multiplication
the {\em convolution product}\/ of reduced functionals.

Assume now that $(\Aa,\del)$ is a regular multiplier Hopf algebra with (non 
trivial) invariant functionals and let $\ph$
be a left invariant functional on $\Aa$. It is known (\cite[Prop.~4.2]{afgd}) 
that the set
\[
\Aah=\bigl\{a\ph:\:a\in\Aa\bigr\}=\bigl\{\ph\,a:\:a\in\Aa\bigr\}
\] 
is also an associative algebra under the convolution product and this product
is non degenerate.

\begin{Prop}
Let $(\Aa,\del)$ be a regular multiplier Hopf algebra with invariant
functionals. Then $\Aah$ embeds into $\Aa^\st$ and $\Aa^\st$ is a subalgebra of 
the multiplier algebra of $\Aah$.
\end{Prop}

\proof
For the first statement we use the fact that $\Aa=\Aa^2$ (this is true for any
multiplier Hopf algebra) and the existence of the modular automorphism for $\ph$, 
i.e.~an automorphism $\sigma$ of $\Aa$ such that
\[
\ph\,a=\sigma(a)\ph
\]
for all $a\in\Aa$ (\cite[Prop.~3.12]{afgd}). Let $b\in\Aa$ and write $b$ as
\[
b=\sum_{k=1}^{N}r_ks_k
\]
with $r_k,s_k\in\Aa$. Then 
\[
\ph\,b=\sum_{k=1}^{N}\sigma(r_k)\ph\,s_k\in\Aa^\st.
\]

Now that we have established that $\Aah\subset\Aa^\st$ we shall use a standard
technique to show that $\Aah$ is in fact an ideal in $\Aa^\st$. Let
$a,b,c\in\Aa$ and $f\in\Aa^\sh$. Set $\gamma=afb$, $\xi=c\ph$.
Further write 
\[
a\tens c=\sum_{k=1}^{N}\del(r_k)(s_k\tens I)
\]
with $r_k,s_k\in\Aa$ for $k=1,\ldots,N$ (this is possible for a regular
multiplier Hopf algebra). Then we have
\[
\begin{array}{r@{\;=\;}l@{\smallskip}}
\bigl(\gamma\conv\xi\bigr)(x)&(afb\tens c\ph)\del(x)\\
&(fb\tens\ph)\bigl(\del(x)(a\tens c)\bigr)\\
&(fb)\left[(\id\tens\ph)\Sum_{k=1}^{N}\del(xr_k)(s_k\tens I)\right]\\
&(fb)\left[\Sum_{k=1}^{N}\bigl[(\id\tens\ph)\del(xr_k)\bigr]s_k\right]\\
&(fb)\left[\Sum_{k=1}^{N}\ph(xr_k)s_k\right]\\
&\Sum_{k=1}^{N}f(bs_k)\bigl(r_k\ph\bigr)(x)
\end{array}
\]
for any $x\in\Aa$, where in the second last equality we used the left invariance
of $\ph$. This way we showed that $\gamma\conv\xi$ belongs to $\Aah$. For the
product $\xi\conv\gamma$ we use the fact that any element of $\Aah$ can be 
written as $d\psi$ with $\psi$ a right invariant functional. Then the argument 
can be repeated with the difference that $d\tens a$ has to be written as
\[
d\tens a=\sum_{k=1}^{N}\del(p_k)(I\tens q_k)
\] 
with $p_k,q_k\in\Aa$, $k=1,\ldots,N$.

Now the associativity of the (convolution) product in $\Aa^\st$ shows that
\[
\xi_1\conv(\gamma\conv\xi_2)=(\xi_1\conv\gamma)\conv\xi_2
\]
for $\gamma\in\Aa^\st$ and $\xi_1,\xi_2\in\Aah$. This means that $\gamma$ is a
multiplier of $\Aah$.
\qed

\section{A characterization of tensor products of multipliers}

Let $(\Aa,\del)$ be a discrete quantum group. Then
the algebra $\Aah$ has a unit (\cite[Prop.~5.3]{afgd}) and consequently 
$\Aah=\Aa^\st$. If in addition $(\Aa,\del)$ is a discrete quantum group, it 
is easy to see that
\[
\Aa^\st=\Span{\bigl.af:\:f\in\Aa^\sh,\:a\in\Aa\bigr.}
=\Span{\bigl.fa:\:f\in\Aa^\sh,\:a\in\Aa\bigr.}.
\] 
In particular for a discrete quantum group $(\Aa,\del)$ we have 
\begin{equation}\label{W}
\Aah=\Span{\bigl.af:\:f\in\Aa^\sh,\:a\in\Aa\bigr.}
=\Span{\bigl.fa:\:f\in\Aa^\sh,\:a\in\Aa\bigr.}.
\end{equation}
Indeed for $a,b\in\Aa$ and $f\in\Aa^\sh$ we have $af=afe_a$ and $fb=e_bfb$
where $e_a$ and $e_b$ are, for example, the units of the ideals generated by
$a$ and $b$ respectively.
Property \refeq{W} is crucial for our next result.   

\begin{Thm}\label{char}
Let $(\Aa,\del)$ be a discrete quantum group and let $\Bb$ be an algebra with 
non degenerate product. Let $Y$ be a multiplier of $\Bb\tens\Aa$. Then
\begin{equation}\label{findim}
\Bigl(\:Y\in\M{\Bb}\tens\M{\Aa}\:\Bigr)\Longleftrightarrow
\Bigl(\:\dim\left\{(\id\tens\xi)(Y):\:\xi\in\Aah\,\right\}<\infty\:\Bigr)
\end{equation}
\end{Thm}

\proof
The implication ``$\Rightarrow$'' is straightforward. Assume that the right
hand side of \refeq{findim} holds. Then let $\{x_1,\ldots,x_N\}$ be a basis in
the space of right slices of $Y$. Fix $b\in\Bb$ and $a\in\Aa$. We shall
consider the expression
\begin{equation}\label{l1}
(\id\tens\ph)\bigl(Y(b\tens a)\bigr)=\bigl((\id\tens a\ph)Y\bigr)b.
\end{equation}
There exist unique scalars $\lambda_k(a)$ such that the right hand side of
\refeq{l1} is equal to
\begin{equation}\label{l2}
\left(\sum_{k=1}^{N}\lambda_k(a)x_k\right)b=\sum_{k=1}^{N}\lambda_k(a)x_kb.
\end{equation}
Clearly each $\lambda_k$ is a linear functional on $\Aa$. By \refeq{W} for any
$c\in\Aa$ we have $\lambda_k c\in\Aah$. By \cite[Lemma 4.11]{afgd} and the 
biduality theorem (\cite[Thm.~4.12]{afgd}) that the map
\[
\Aah\ni d\ph\longmapsto\bigl(\lambda_kc\bigr)(d)=\lambda_k(cd)
\]
is determined by a unique element $q\in\Aa$ in such a way that
\begin{equation}\label{osta}
\lambda_k(cd)=\bigl(d\ph\bigr)(q).
\end{equation}
Writing $q=y_k(c)$ we define linear maps $y_k\colon\Aa\to\Aa$ for
$k=1,\ldots,N$.

Consider now
\[
(\id\tens\ph)\bigl((b\tens a)Y)=b\bigl((\id\tens\ph\,a)Y\bigr)=
b\bigl((\id\tens\sigma(a)\ph)Y\bigr).
\]
As before we write this last expression as
\[
b\left(\sum_{k=1}^{N}\lambda_k\bigl(\sigma(a)\bigr)x_k\right)=
\sum_{k=1}^{N}\lambda_k\bigl(\sigma(a)\bigr)bx_k
\]
and we can define linear maps $z_k\colon\Aa\to\Aa$ by
\[
\bigl(\ph\,d\bigr)\bigl(z_k(c)\bigr)=\lambda_k\bigl(\sigma(dc)\bigr)
\]
or in other words
\[
\bigl(d\ph\bigr)\bigl(z_k(c)\bigr)=\lambda_k\bigl(d\sigma(c)\bigr).
\]
We shall show that 
$\bigl((y_k,z_k)\bigr)_{k=1,\ldots,N}$ are double centralizers of
$\Aa$. Take $a_1,a_2\in\Aa$ and $d\ph\in\Aah$. We have
\[
\begin{array}{r@{\;=\;}l@{\smallskip}}
\bigl(d\ph\bigr)\bigl(z_k(a_1)a_2\bigr)
&\bigl((a_2d)\ph\bigr)\bigl(z_k(a_1)\bigr)
=\lambda_k\bigl(a_2d\sigma(a_1)\bigr)\\
\bigl(d\ph\bigr)\bigl(a_1y_k(a_2)\bigr)
&\ph\bigl(a_1y_k(a_2)d\bigr)=\ph\bigl(y_k(a_2)d\sigma(a_1)\bigr)\\
&\left[\bigl(d\sigma(a_1)\bigr)\ph\right]\bigl(y_k(a_2)\bigr)
=\lambda_k\bigl(a_2d\sigma(a_1)\bigr)
\end{array}  
\]
and since this equality holds for any $d\ph\in\Aah$ we have that
$z_k(a_1)a_2=a_1y_k(a_2)$ for all $a_1,a_2\in\Aa$. This way we have obtained
multipliers $(y_k)_{k=1,\ldots,N}$ of $\Aa$. 

Upon substitution of $y_kc$ for $q$ in equation \refeq{osta} we obtain
\[
\lambda_k(cd)=\ph(y_kcd)
\]
for all $c,d\in\Aa$. Inserting $c=a$ and $d$ equal to the unit of the ideal
generated by $a$ we obtain $\lambda_k(a)=\ph(y_ka)$ for $k=1,\ldots,N$.
Therefore (cf.~\refeq{l1} and \refeq{l2}) the multipliers 
$(y_k)_{k=1,\ldots,N}$ satisfy
\[
(\id\tens\ph)\left[\left(\sum_{k=1}^{N}x_k\tens y_k\right)(b\tens a)\right]
=(\id\tens\ph)\bigl[Y(b\tens a)\bigr].
\]
As $a\in\Aa$ and $b\in\Bb$ were arbitrary we obtain this equality 
for all $a\in\Aa$ and $b\in\Bb$. By faithfulness of $\ph$
(cf.~\cite[Sect.~3]{afgd}) we have
\[
Y=\sum_{k=1}^{N}x_k\tens y_k
\]
and consequently $Y\in\M{\Bb}\tens\M{\Aa}$.
\qed

In \cite{ap} the following definition of an almost periodic element for a 
discrete quantum group $(\Aa,\del)$ was proposed: a multiplier $x\in\M{\Aa}$
is an almost periodic element for $(\Aa,\del)$ if
$\del(x)\in\M{\Aa}\tens\M{\Aa}$ (where we use the unique extension of $\del$
to the multiplier algebra). It was shown that the set of all almost periodic
elements for $(\Aa,\delta)$ forms a Hopf $*$-algebra. On the formal level the
definition of an almost periodic element is similar to that of an almost
periodic function on a locally compact group (cf.~\cite[\S 41]{loom}).   
Theorem \ref{char} shows that this definition corresponds more directly to the
notion of an {\em almost invariant} function (\cite[Sect.~39D]{loom}), but
this distinction fades on the purely algebraic level.

\section{Acknowledgements}
The author is very grateful to professor Alfons Van Daele for invaluable help
and encouragement. Also he would like to extend his acknowledgements to
professor Joachim Cuntz and colleagues from the Mathematisches Institut of the 
University of M\"unster for creating an outstanding scientific environment.


\begin{thebibliography}{XX}
\bibitem{loom}
{\sc Loomis, L.H.} {\em An introduction to abstract harmonic analysis,}\/ Van
Nostrad, Princeton N.J.~1953.
%\bibitem{mnw}
%{\sc Masuda, T., Nakagami, Y.~\& Woronowicz, S.L.} A $C^*$-algebraic framework
%for quantum groups, to appear in {\em Int.~J.~Math.}
\bibitem{pw}
{\sc Podle\'{s}, P.~\& Woronowicz, S.L.} Quantum deformation of Lorentz group,
{\em Comm. Math.~Phys.}\/ \textbf{130} (1990), 381--431.
\bibitem{ap}
{\sc So\l{}tan, P.M.} Compactifications of discrete quantum groups, Preprint
des SFB 478 des Mathematisches Institut der Wesf\"alischen
Wilhelms-Universit\"at M\"unster (2003).
%\bibitem{dp}
%{\sc Van Daele, A.} Dual pairs of Hopf $*$-algebras, {\em Bull.~London
%Math.~Soc.}\/ \textbf{25} (1993), 209--230.
\bibitem{mha}
{\sc Van Daele, A.} Multiplier Hopf algebras, {\em Trans.~AMS}\/ \textbf{342}
No.~2 (1994), 917--932.
\bibitem{dqg}
{\sc Van Daele, A.} Discrete quantum groups, {\em J.~Algebra}\/ \textbf{180} 
(1996), 431--444.
\bibitem{afgd}
{\sc Van Daele, A.} An algebraic framework for group duality, {\em Adv.~in
Math.}\/ \textbf{140} (1998), 323--366.
%\bibitem{dt}
%{\sc Van Daele, A.~\& Zhang, Y.} Multiplier Hopf algebras of discrete type,
%{\em J.~Algebra}\/ \textbf{214} (1999), 400--417.
\end{thebibliography}
\end{document}